\documentclass[a4paper,12pt]{amsart}

\usepackage[utf8]{inputenc}
\usepackage{amsmath, amsfonts, amssymb, amsthm, bm}
\usepackage[dvipsnames]{xcolor}
\usepackage{xcolor}
\usepackage{mathrsfs}
\usepackage{amsfonts}
\usepackage{amsmath, amsfonts, amssymb, amsthm, amscd}
\usepackage{lmodern}
\usepackage{lmodern}
\usepackage{amsthm}

\newcommand{\SL}{\operatorname{SL}}

\newtheorem{theorem}{Theorem}[]

\newtheorem{remark}[theorem]{Remark}

\newtheorem{corollary}[theorem]{Corollary}
\newtheorem{lemma}[theorem]{Lemma}

\newtheorem{acknowledgement}{Acknowledgement}

\def\R{\mathbb R} \def\Z{\mathbb Z}

\def\Q{\mathbb Q}

\def\F{\mathbb F}

\def\build#1_#2^#3{\mathrel{\mathop{\kern 0pt#1}\limits_{#2}^{#3}}}

\usepackage{color}
\usepackage{hyperref}
\usepackage{letltxmacro}
\hypersetup{
    colorlinks=true, 
    linktoc=all, 
    linkcolor=blue, 
   citecolor=blue,
    filecolor=blue,
    urlcolor=blue
}
\usepackage[mathlines]{lineno} 
\usepackage{amsmath}

\usepackage{etoolbox} 

\title{On values of isotropic quadratic forms}
\author{Manoj Choudhuri and Prashant J. Makadiya}
\address{Institute of Infrastructure, Technology, Research and Management, Near Khokhara Circle, maninagar (East), Ahmedabad 380026, Gujarat, India.}
\email{manojchoudhuri@iitram.ac.in} \email{prashant.makadiya.20pm@iitram.ac.in}
\date{}

\begin{document}

\begin{abstract}
Let $K$ be a locally compact non-discrete field of characteristic $p>2$ and $Q$ be a non-degenerate isotropic binary quadratic form with coefficients in $K$. We obtain asymptotic estimates for the number of solutions in the two-fold product of a discrete subring inside $K$, of the inequalities of the form $|Q(x,y)|<\delta$ for some $\delta>0$, where $| \cdot |$ is an ultrametric absolute value on $K$. The estimates are obtained in terms of continued fraction expansions of the coefficients of the quadratic form $Q$.   
\end{abstract}

\maketitle

\smallskip

{\it Mathematics Subject Classification,} Primary: $11$E$08$; Secondary: $11$J$61$, $11$J$70$, $11$J$83$
    $11$K$50$, $37$A$44$.\\
{\it Keywords:} Quadratic forms, locally compact fields, asymptotic estimates, continued fractions.
\maketitle


\vspace{0.2cm}

The Oppenheim conjecture, solved by Margulis in 1987 (see \cite{M} for more details), states that if $Q$ is a real non-degenerate indefinite quadratic form which is not proportional to a form with rational coefficients, then $Q(\Z^n)$ is dense in $\R$ if $n\geq 3$. After Oppenheim conjecture was settled, people got interested in studying finer questions related to the distribution of the values of $Q$ on integral points. Given a quadratic form as above, and $a,b,\rho\in\R$ with $\rho>0$, let 
\begin{align*}
N_Q(a,b,\rho):=\#\ \{v\in\Z^n:\ a<Q(v)<b,\ v\in B(\rho)\},
\end{align*}
$B(\rho)$ being the ball of radius $\rho$ around the origin in $\R^n$. Also let
\begin{align*}
V_Q(a,b,\rho):=\text{Vol}\ (\{v\in\R^n:\ a<Q(v)<b,\ v\in B(\rho)\}).
\end{align*}
Then it was shown by Dani and Margulis in \cite{DM1} that 
\begin{align*}
\liminf\limits_{\rho\rightarrow\infty}\frac{N_Q(a,b,\rho)}{V_Q(a,b,\rho)}=1.    
\end{align*}
Asymptotic upper bound for the quantity $\frac{N_Q(a,b,\rho)}{V_Q(a,b,\rho)}$ was found by Eskin, Margulis and Mozes (see \cite{EMM1} for instance), and combining the result of \cite{DM1}, they showed that if $Q$ is a quadratic form as above such that the signature of $Q$ is neither $(2,1)$ nor $(2,2)$, then 
\begin{align*}
\lim_{\rho\rightarrow\infty}\frac{N_Q(a,b,\rho)}{V_Q(a,b,\rho)}=1.    
\end{align*}

The Oppenheim conjecture fails for binary quadratic forms due to the existence of badly approximable numbers. A real number $\alpha$ is called badly approximable if there exists $c>0$ such that $\left|\alpha-\frac{\displaystyle p}{\displaystyle q}\right|>\frac{\displaystyle{c}}{\displaystyle{q^2}}$ for any rational number $\frac{\displaystyle p}{\displaystyle q}$. Now, let $Q$ be the binary quadratic form defined by 
\begin{align}\label{bad}
Q(x,y)=(x+\alpha y)y,  
\end{align}
$\alpha$ being a badly approximable number. Then $Q(\Z^2)$ avoids the neighbourhood $(-c,c)$ of zero. Nevertheless, one can study the distribution of the values taken by such forms at integral points. This was done in \cite{CD1} with the interval $(a,b)$ being a neighbourhood of $0$. In case of binary quadratic forms, the asymptotic estimates depend on the quadratic form under consideration, and they are given in terms of the partial quotients of the continued fraction expansions of the coefficients of the quadratic form. There is a natural connection between the values of non-degenerate indefinite binary quadratic forms at integral points, and certain geometric and dynamical aspects of the orbits of geodesic flow associated with the modular surface. In \cite{CD1}, the authors explored this connection, and used a method of coding of geodesics on the modular surface via nearest integer continued fraction which was introduced by S. Katok and I. Ugarcovicci (see \cite{KU} for instance), to obtain the estimates (see \cite{Si1} for a different proof which does not use the mechinary of geodesic flow etc.). The method of \cite{CD1} can be adopted to obtain similar type of estimates in terms of a more general class of continued farctions as well, see Remark $3.4$ of \cite{C1} for more details. 

In the present article, we confine our attention to binary isotropic quadratic forms with coefficients in a non-discrete locally compact field $K$ of characteristic greater than $2$. Without loss of generality, we may assume that $K$ is the formal Laurent series field in one indeterminant over a finite field (see Theorem $8$ in chapter $I$ of \cite{WA} for instance). In other words, we may assume that  $K:=\ \F_q((X^{-1}))$, where $\F_q$ is the finite field with $q$ elements with $q=p^r$ for some prime $p>2$ and positive integer $r$. In this set up, the role of $\Z$ is played by the polynomial ring $\F_q[X]$, and the role of $\Q$ is played by $\F_q(X)$. For simplicity, we will use the notation $Z$ to denote the polynomial ring $\F_q[X]$. 

One can define a valuation $\nu$ on $K$ as follows: if $\alpha=\sum\limits_{j\geq j_0 }a_jX^{-j}\in K$ with $a_j\in\F_q$, then 
\begin{align*}
\nu(\alpha):=\ \text{inf}\ \{j\in\Z:a_j\neq 0\}.  
\end{align*}
This valuation gives rise to an absolute value on $K$ as follows: if $\alpha(\neq 0)\in K$ and $\nu(\alpha)=d_{\alpha}$, then
\begin{align*}
|\alpha|:=q^{d_{\alpha}},
\end{align*}
and the absolute value of the zero element in $K$ is $0$. Then $K$ is the completion of $\F_q(X)$ with respect to this absolute value. As $\nu$ is a non-Archimedean valuation, the absolute value defined above is an ultrametric absolute value. In the metric generated by this absolute value, $Z=\F_q[X]$ is a discrete subring of $K$.

We call an element $\alpha\in K$ irrational if $\alpha\notin \F_q(X)$. The notion of badly approximable elements can also be defined as in the case of real numbers. In particular, an element $\alpha\in K$ is said to be badly approximable if there exista a real constant $c>0$ such that $$\left|\alpha - \frac{\displaystyle P}{\displaystyle Q}\right|>\frac{\displaystyle c}{\displaystyle{|Q|^2}}$$ for all $P,Q\in Z$. Then considering the same example as in (\ref{bad}) with $\alpha$ a badly approximable element, it is easy to see that Oppenheim conjecture fails for binary quadratic forms with coefficients in $K$ as well. In this article, we do a similar study as in \cite{CD1} for non-denerate isotropic binary quadratic forms with coefficients in $K$. 

Being a locally compact field, $K$ admits a Haar measure (see \cite{RV1} for details) which we denote by $\mu$. For $a\in K$ and $r\in\Z$, let
\begin{align*}
B(a,q^r):=\{\alpha\in K\ :\ |\alpha-a| < q^r\}
\end{align*}
be the open disc around $a$ of radius $q^r$, then $\mu(B(a,q^r))=q^r$. Also, let $\mu\otimes\mu$ be the corresponding product measure on $K^2$ which is denoted by $\eta$.

As in the case of real numbers, any $\alpha$ in $K$ has a unique continued fraction expansion 
\begin{align*}
 \alpha=b_0+\frac{1}{\displaystyle{b_1+\frac{1}{\displaystyle{b_2+\frac{1}{\displaystyle{b_3+....}}}}}},   
\end{align*}
also written as
\begin{align*}
\alpha=[b_0,b_1,b_2,....]
\end{align*}
with $b_j\in Z$ for $j\geq 0$ and $b_j$ has positive degree for $j\geq 1$. Given any $\alpha=\displaystyle{\sum\limits_{j\geq j_0}a_jX^{-j}}$ in $K$, let
\begin{align*}
\lfloor \alpha \rfloor =\displaystyle{\left\{ \begin{array}{rcl} \displaystyle{\sum \limits_{j=j_0}^{0}a_{j}X^{-j}} & \mbox{if} & j_0 \leq 0  \\ 0 \ \ \ & \mbox{if} & j_0 \geq 1.
\end{array}\right.}
\end{align*}
Then the continued fraction algorithm is defined as follows:
\begin{align*}
\alpha_0:=\alpha,\ \alpha_{n+1}:=(\alpha_n -b_n)^{-1}\ \text{and}\ b_n=\lfloor \alpha_n \rfloor.    
\end{align*}
Here, $b_n$'s are called partial quotients and $\alpha_n$'s are called complete quotients of the continued fraction expansion of $\alpha$ (see \cite{WSch1} for more details).
Now, let $\frac{\displaystyle{s_n}}{\displaystyle{t_n}}$ be the $n$th convergent of the continued fraction expansion of $\alpha$, i.e.,
\begin{align*}
\frac{s_n}{t_n}=[b_0,b_1,b_2,...,b_n].
\end{align*}
Then the sequences $(s_n)_{n\geq 0}$ and $(t_n)_{n\geq 0}$ in $Z$ satisfy the following recurrence relations: 
\begin{equation}\label{eq1}
 s_n=b_ns_{n-1}+s_{n-2},\ \ t_n=b_nt_{n-1}+t_{n-2}. 
\end{equation}
They also satisfy the following equation:
\begin{equation}\label{eq2}
s_{n+1}t_n-s_nt_{n+1}=(-1)^n 
\end{equation}
which tells us that $s_n$ and $t_n$ are coprime, i.e., they do not have any common factor other than the constant polynomials in $\F_q[X]$.
The following equalities which are special features of continued fraction theory, will be quite useful for this article. If $\alpha,\ b_n,\  s_n,\ t_n$ are as above, then 
\begin{equation}\label{eq4}
|t_{n}|=|b_{n} \cdots b_{1}| \ ; \forall n \geq 1,
\end{equation}

\begin{equation}\label{eq3}
 \left|\alpha - \frac{s_n}{t_n}\right|\ =\ \frac{1}{|b_{n+1}||t_n|^2},
 \end{equation}
 and 
 \begin{equation}\label{eqn}
 \left|\alpha - \frac{s_n}{t_n}\right|\ =\ \frac{1}{|t_{n+1}||t_n|}.    
 \end{equation}
Note that in the case of continued fraction for real numbers, inequalities hold instead of equalities in  $(\ref{eq3})$ and $(\ref{eqn})$. This is because of the ultrametric nature of the absolute value on $K$.
The following lemma is a simple characterization of the convergents of the continued fraction expansion of any element in $K$, the proof of which can be found in \cite{WSch1}.
\begin{lemma}\label{lem1}
 Let $s,t\in Z$ with $t\neq 0$. Then $\frac{\displaystyle{s}}{\displaystyle{t}}$ is a convergent to $\alpha$ if and only if  
 \begin{equation}\label{eq5}
 \left|\alpha - \frac{s}{t}\right| < \frac{1}{|t|^2}.
 \end{equation}
\end{lemma}
Now, let us consider binary quadratic forms with coefficients in $K$.
It is well-known that if $Q$ is a non-degenerate isotropic quadratic form with coefficients in a field $F$ of characteristic not equal to $2$, then there exists a basis $\{v_1,v_2\}$ of $F^2$ such that if $a_1,a_2\in F$, then 
\begin{align*}
Q(a_1v_1+a_2v_2)=a_1a_2.
\end{align*}
This says in particular that if $Q_0$ is the quadratic from on $K^2$ defined by
\begin{align*}
Q_0(x,y)=xy\ \text{for}\ x,y\in K,
\end{align*}
then for any non-degenerate isotropic quadratic form $Q$ on $K^2$, there is a matrix $A_Q$ in $ \SL(2,K)$ and $\gamma$ in $K$, such that
\begin{equation}\label{eqnew}
Q(x,y)=\gamma\ Q_0(A_Q(x,y)).
\end{equation}
So, to study the asymptotic behaviour of the set of values of an isotropic quadratic form with coefficients in $K,$ it is enough to consider quadratic form $Q$ given as follows:
\begin{align*}
Q(x,y)=(ax+by)(cx+dy)  
\end{align*}
with $a,b,c,d \in K, \ bc-ad=1$ (there is no loss of generality because one may replace $\gamma$ by $-\gamma$ in (\ref{eqnew})).

Now, let $Q$ be a quadratic form of the type $Q(x,y)=(ax+by)(cx+dy)$ with $a,b,c,d \in K \  \text{and}\ bc-ad=1$,  such that $\frac{\displaystyle{b}}{\displaystyle{a}}$ is an irrational element of $K$.
Also, let $\mathfrak{p}$ be the set of primitive elements of $Z^{2}$, i.e.,
$\mathfrak{p}$ is the set of those $(s,t)$ in $Z^2$ such that $s$ and $t$ do not have a common factor except constant polynomials.
For fixed real numbers $k$ and $\delta$ with $k>1$ and $0<\delta <1$, let
\begin{align*}
G(\rho):= \{ (s,t) \in \mathfrak{p} : 0<|Q(s,t)|<\delta, \ ||(s,t)|| \leq \rho, \ |cs+dt| >k \},
\end{align*}
where $||(s,t)||=\max\ \{ |s|,|t| \}$. Let $\alpha=- \frac{\displaystyle b}{\displaystyle a}$ and $\beta =ac$, and the continued fraction expansion of $\alpha$ be given by
\begin{align*}
\alpha=[b_{0},b_{1},b_{2},...]
\end{align*}
with $\frac{\displaystyle{s_n}}{\displaystyle{t_n}}$ being the $n$th convergent. Also let 
\begin{align*}
H(\rho):= \{ (x,y) \in K^{2} : 0<|Q(x,y)|<\delta, \ ||(x,y)|| \leq \rho, \ |cx+dy| >k \}.    
\end{align*}
In this article, we find asymptotic lower and upper bound of the quotient  $\frac{\displaystyle{\#\  G(\rho)}}{\displaystyle{\eta\ (H(\rho))}}$ as $\rho \rightarrow \infty$. Let
\begin{align*}
\alpha^-:=\liminf \limits_{n\rightarrow\infty}\ \frac{1}{n}\sum\limits_{j=1}^{n}\deg\left(b_j\right)
\end{align*}
and 
\begin{align*}
\alpha^+:=\limsup \limits_{n\rightarrow\infty}\ \frac{1}{n}\sum\limits_{j=1}^{n}\deg \left(b_j\right),     
\end{align*} 
where $\text{deg}\left(b_j\right)$ denotes the degree of the polynomial $b_j$.
Also for $0<\delta<1$, let 
\begin{align*}
e(\delta):=\liminf \limits_{n\rightarrow\infty}\ \frac{1}{n}\#\left\{j,\ 1\leq j\leq n:\ |b_{j+1}|\geq \frac{1}{\delta}\right\}    
\end{align*} 
and 
\begin{align*}
f(\delta):=\limsup \limits_{n\rightarrow\infty}\ \frac{1}{n}\#\left\{j,\ 1\leq j\leq n:\ |b_{j+1}|\geq \frac{1}{\delta}\right\}.   
\end{align*} 
Then the main result of this article is contained in the following theorem.
\begin{theorem}\label{th2}
Let $Q$ be a quadratic form defined by 
\begin{align*}
Q(x,y)=(ax+by)(cx+dy)\ \text{with}\  a,b,c,d\in K,\ bc-ad=1,   
\end{align*} 
and $\frac{\displaystyle b}{\displaystyle a}$ an irrational element of $K$. Also let $G(\rho)$, $H(\rho)$, $\alpha^{+}$, $\alpha^{-}$, $e(\delta)$, $f(\delta)$ be as defined above. If $\alpha^-<\infty$, then we have the followings: 
\begin{align*}
 \liminf\limits_{\rho\rightarrow\infty}\frac{\displaystyle{\#\ G(\rho)}}{\displaystyle{\eta\ (H(\rho))}}\ \geq\ \mathfrak{c}\ \frac{\displaystyle{e(\delta)}}{\displaystyle{\alpha^+}} 
\end{align*} 
and
\begin{align*}
\limsup\limits_{\rho\rightarrow\infty}\frac{\displaystyle{\#\ G(\rho)}}{\displaystyle{\eta\ (H(\rho))}}\ \leq\ \mathfrak{c}\  \frac{\displaystyle{f(\delta)}}{\displaystyle{\alpha^-}},  
\end{align*}
 where $\mathfrak{c}$ is a constant depending on $\delta$ and $q$.
\end{theorem}

\begin{remark}\label{rem3}
Let 
\begin{align*}
I(\rho):= \{ (s,t) \in \mathfrak{p} \ : \ 0<|Q(s,t)|< \delta, \ ||(s,t)|| \leq \rho, \ |as+bt|>k \}   
\end{align*}
and
\begin{align*}
 J(\rho):= \{ (x,y) \in K^{2} \ : \ 0<|Q(x,y)|< \delta, \ ||(x,y)|| \leq \rho, \ |ax+by|>k \}.   
\end{align*}
Then one can obtain a similar estimates for $\frac{\displaystyle{\#\ I(\rho)}}{\displaystyle{\eta\ (J(\rho))}}$
in terms of the continued fraction expansion of $- \frac{\displaystyle d}{\displaystyle c}$ provided $\frac{\displaystyle d}{\displaystyle c}$ is an irrational element of $K$.
\end{remark}
\noindent Proof of Theorem \ref{th2}:\\
Let 
\begin{align*}
G'(\rho):=\{ (s,t) \in \mathfrak{p} : |t(t \alpha - s)|< \delta, \ |t|\leq \rho \}.  
\end{align*}
It is easy to see that
\begin{equation}\label{eq7}
Q(s,t)=(t \alpha - s)(t+ \beta (t \alpha - s)).
\end{equation}
If $|Q(s,t)|<\delta$ with $|cs+dt|>k$ then $|as+bt|<\frac{\displaystyle{\delta}}{\displaystyle k}$, which implies that $|t \alpha - s|< \frac{\displaystyle{\delta |a|}}{\displaystyle k}$, i.e., $|t \alpha - s|$ is bounded. Now by (\ref{eq7}),
\begin{align*}
\frac{\displaystyle{|Q(s,t)|}}{\displaystyle{|t(t \alpha - s)|}}=\left| 1+\frac{\displaystyle{\beta}}{\displaystyle t}(t \alpha - s) \right|.   
\end{align*} 
Since $|t \alpha - s|$ is bounded, it follows that $\frac{\displaystyle{|Q(s,t)|}}{\displaystyle{|t(t \alpha - s)|}}=1$ if $|t|$ is sufficiently large. Note that when $|t \alpha - s|$ is bounded, $||(s,t)|| \rightarrow \infty$ if and only if $|t| \rightarrow \infty$. Also, if $|t(t \alpha - s)|< \delta$, then clearly $|t \alpha - s|$ is bounded and $\frac{\displaystyle{|Q(s,t)|}}{\displaystyle{|t(t \alpha - s)|}}=1$ for sufficiently large $|t|$. Combining all these facts, we can say that there exists a constant $C>0$ such that 
\begin{align*}
\# G^{'}(\rho)-C \leq \# G(\rho) \leq \# G^{'}(\rho)+C  
\end{align*}
for sufficiently large $\rho$. Since $0<\delta <1$, it follows from Lemma \ref{lem1}, that if $(s,t) \in G^{'}(\rho)$, then $s=s_{j}$ and $t=t_{j}$, where $\frac{\displaystyle{s_{j}}}{\displaystyle{t_{j}}}$ is a convergent of $\alpha$ in its continued fraction expansion. Also $G^{'}(\rho)=G^{'}(|t_{n}|)$ if $|t_{n}| \leq \rho <|t_{n+1}|$. Note that if $(s_{j},t_{j}) \in G^{'}(|t_{n}|)$, then $(as_{j},at_{j}) \in G^{'}(|t_{n}|)$ as well for any $a \in \F^{*}_q$. 

Now, let us calculate the measure of $H(\rho)$. Let $A$ be the set given by 
\begin{align*}
A:=\{ (x,y)\in K^{2}:0<|xy|<\delta, \ ||(x,y)|| \leq \rho , \ |y|>k \},    
\end{align*}
then 
\begin{align*}
\eta(H(\rho))=|\text{det}(M)|\ \eta(A)  
\end{align*}
where $M=\left[ \begin{array}{crcr} a & b \\ c & d \end{array} \right]$. Since $bc-ad=1$, we have that $\eta(H(\rho))=\eta(A)$.

Note that for $0<\delta <1, \ k>1$ and $\rho \geq k$, there exist unique $m_{0}, \ m^{'}_{0}, \ t$ and $i \in \mathbb{Z}$ such that $ q^{m_{0}}\leq \delta < q^{m_{0}+1}, \ q^{m^{'}_{0}}\leq \sqrt{ \delta } < q^{m^{'}_{0}+1}, \ q^{m^{'}_{0}+t}\leq k < q^{m^{'}_{0}+t+1}$ and $q^{m^{'}_{0}+t+i}\leq \rho < q^{m^{'}_{0}+t+i+1}$. Also for $1 \leq n \leq i$, let 
\begin{align*}
A_{n}:=\{ (x,y) \in K^{2}: |x| \leq q^{m_{0}-m^{'}_{0}-t-n} \ and \ |y|=q^{m^{'}_{0}+t+n} \}.
\end{align*}
Clearly $A_n$'s are disjoint, and it is easy to see that $A=\cup_{n=1}^{i}A_{n}$.
Hence, $\eta (A)=\sum \limits_{n=1}^{i}\eta (A_{n})$. Now, 
\begin{align*}
\{ y \in  K : |y| \leq q^{m^{'}_{0}+t+n} \}\\& =\{ y \in  K : |y| < q^{m^{'}_{0}+t+n} \} \cup \{ y \in  K : |y| = q^{m^{'}_{0}+t+n} \}.
\end{align*}
Therefore,
\begin{align*}
\eta (A_{n})&=\mu ( \{ x \in  K : |x| \leq q^{m_{0}-m^{'}_{0}-t-n} \} ) \cdot \mu ( \{ y \in  K : |y|= q^{m^{'}_{0}+t+n}  \} )\\&=\mu ( \{ x \in  K : |x| \leq q^{m_{0}-m^{'}_{0}-t-n} \}) \\& \ \ \ \cdot (\mu ( \{ y \in  K : |y| \leq q^{m^{'}_{0}+t+n} \}) - \mu ( \{ y \in  K : |y| < q^{m^{'}_{0}+t+n} \} ))\\& =(q^{m_{0}-m^{'}_{0}-t-n+1}) \cdot (q^{m^{'}_{0}+t+n+1}-q^{m^{'}_{0}+t+n})\\&=(q^{m_{0}-m^{'}_{0}-t-n+1})(q^{m^{'}_{0}+t+n})(q-1)\\&=q^{m_{0}+1}(q-1),
\end{align*}
and consequently,
\begin{align*}
\eta (H(\rho))=\eta (A)=\sum\limits_{n=1}^{i}\eta(A_{n})=iq^{m_{0}+1}(q-1).
\end{align*}
As $q^{m'_0+t+i}\leq\rho < q^{m'_0+t+i+1}$, taking logarithm with base $q$ we obtain
\begin{align*}
(m'_0+t+i)\leq\log_q\rho < (m'_0+t+i+1).
\end{align*}
Consequently,
\begin{align*}
\log_q\rho-m'_0-t-1 <i \leq \log_q \rho-m'_0-t.
\end{align*}
Then it follows that
\begin{align}\label{eq8}
\left( \log_q \rho-m^{'}_{0}-t-1 \right)(q-1)q^{m_{0}+1}\\  &<& \eta (H(\rho))  \leq   \left( \log_q \rho-m^{'}_{0}-t \right)(q-1)q^{m_{0}+1} \nonumber .
\end{align}
Now,
\begin{align*}
\liminf \limits_{\rho \rightarrow \infty} \frac{\displaystyle{\# G(\rho)}}{\displaystyle{\eta (H(\rho))}}& \geq \liminf \limits_{\rho \rightarrow \infty} \frac{\displaystyle{\# G'(\rho)-C}}{\displaystyle{\eta (H(\rho))}}\\ & =\liminf \limits_{n \rightarrow \infty} \frac{\displaystyle{\# G'(|t_{n}|)-C}}{\displaystyle{\eta (H(|t_{n}|))}}\ \ (\text{for}\ |t_n|\leq \rho< |t_{n+1}|)\\& =\liminf \limits_{n \rightarrow \infty} \frac{\displaystyle{\frac{\displaystyle 1}{\displaystyle n} ( \# G'(|t_{n}|)-C)}}{\displaystyle{\frac{\displaystyle 1}{\displaystyle n} ( \eta (H(|t_{n}|)))}}\\& \geq \frac{\liminf \limits_{n \rightarrow \infty}\frac{\displaystyle 1}{\displaystyle n} ( \# G'(|t_{n}|))}{\limsup \limits_{n \rightarrow \infty}\frac{\displaystyle 1}{\displaystyle n} ( \eta (H(|t_{n}|)))}\\& \geq \frac{\displaystyle{\liminf \limits_{n \rightarrow \infty}\frac{\displaystyle 1}{\displaystyle n} ( q-1) \ \# \left\{ j: 1 \leq j \leq n, \ |b_{j}| \geq \frac{1}{\delta} \right\} }} {\displaystyle{\limsup \limits_{n \rightarrow \infty}\frac{\displaystyle 1}{\displaystyle n} \left( \log_q |t_n|-m'_0-t \right) q^{m_0+1}(q-1)}}\  (\text{by}\ (\ref{eq3})\ \text{and}\ (\ref{eq8}))\\& \geq \frac{\displaystyle{\liminf \limits_{n \rightarrow \infty}\frac{\displaystyle 1}{\displaystyle n} \ \# \left\{ j: 1 \leq j \leq n, \ |b_{j}| \geq \frac{1}{\delta} \right\} }} {\displaystyle{\limsup \limits_{n \rightarrow \infty}\frac{\displaystyle 1}{\displaystyle n} \left(\log_q |b_1b_2 \cdots b_n|-m'_{0}-t \right) q^{m_0+1}}}\  (\text{by}\ (\ref{eq4}))\\& \geq \frac{\displaystyle{\liminf \limits_{n \rightarrow \infty}\frac{\displaystyle 1}{\displaystyle n} \ \# \left\{ j: 1 \leq j \leq n, \ |b_{j}| \geq \frac{1}{\delta} \right\} }} {\displaystyle{\limsup \limits_{n \rightarrow \infty}\frac{\displaystyle 1}{\displaystyle n} \left(\sum \limits_{j=1}^{n}\deg(b_j)-m'_0-t \right) q^{m_0+1}}}\\& =\frac{e(\delta)}{\alpha^{+}} q^{-(m_0+1)}.
\end{align*}

A similar calculation yields 
\begin{align*}
\limsup \limits_{\rho \rightarrow \infty} \frac{\# G(\rho)}{\eta (H(\rho))} \leq \frac{f(\delta)}{\alpha^{-}} q^{-(m_0+1)}.
\end{align*}
This completes the proof of the theorem with the value of the constant $\mathfrak{c}$ being $q^{-(m_0+1)}$.
\begin{corollary}\label{cor4}
 Let $Q$ be a quadratic form as in Theorem \ref{th2}, and $0<\delta<1$ be fixed. Then there exist a subset $K'$ of $K$ with $\mu (K')=\mu (K)$ such that if $\alpha =-\frac{\displaystyle b}{\displaystyle a} \in K'$, then
 \begin{align*}
\lim_{\rho\rightarrow\infty}\frac{\# G(\rho)}{\eta (H(\rho))}=\frac{q-1}{{q^{\lceil \delta^{-1} \rceil  +m_0+1 }}},    
\end{align*}
 where $\lceil \delta^{ - 1} \rceil$ denotes the smallest integer greater or equal to $\delta^{ - 1}$. 
\end{corollary}
\begin{proof}
Let $[b_0,b_1,b_2,\dots]$ be the continued fraction expansion of $\alpha = -\frac{\displaystyle{b}}{\displaystyle{a}}$ as above. By Theorem $6$ of \cite{BN1}, there is a full measure subset $K'$ of $K$ such that if $\alpha = -\frac{\displaystyle{b}}{\displaystyle{a}}\in K'$, then 
\begin{align*}
\lim \limits _{n \rightarrow \infty}\frac{1}{n}\sum \limits_{j=1}^{n} \deg\left(b_j\right)={\frac{q}{q-1}},
\end{align*} 
and, therefore, $\alpha^{-}=\alpha^{+}={\frac{q}{q-1}}$. Also for any $0< \delta <1$, there exists a unique $l \in \mathbb{N}$ such that $l =\lceil \delta^{-1} \rceil$. Then by Theorem $14$ of \cite{LN1}, for $\alpha$ in a full measure set which without loss of generality we may assume to be $K'$, 
\begin{align*}
\lim \limits _{n \rightarrow \infty}\frac{1}{n}\ \# \{ 1 \leqslant j \leqslant n : |b_j| \geqslant q^l \} = \frac{\displaystyle 1}{\displaystyle{q^{l-1}}},
\end{align*}
and, therefore, $e(\delta)=f(\delta)= \frac{\displaystyle 1}{\displaystyle{q^{l-1}}}=\frac{\displaystyle 1}{\displaystyle{q^{\lceil \delta^{-1} \rceil -1}}}$. 
Then it follows from Theorem \ref{th2} above that, if $\alpha=-\frac{\displaystyle b}{\displaystyle a}\in K^{'}$, then
\begin{align*} 
\lim_{\rho \rightarrow \infty} \frac{\# G(\rho)}{\eta (H(\rho))}
=\frac{\displaystyle{{q-1}}}{{{q}^{\displaystyle{{\lceil \delta^{-1} \rceil +m_{0}+1}}}}}. 
\end{align*}
\end{proof}
\begin{remark}\label{rem5}
Let $Q$, $\alpha$ be as in Theorem \ref{th2}. Now, if the absolute values of the partial quotients in the continued fraction expansion of $\alpha$ are bounded by some real numbers (equivalently, $\alpha$ is badly approximable, cf. \cite{L1}), then it is easy to see that $e(\delta)=f(\delta)=0$ if $\delta$ is sufficiently small. In this case, 
\begin{align*}
\lim \limits_{\rho \rightarrow \infty} \frac{\# G(\rho)}{\eta (H(\rho))}=0.    
\end{align*}
\end{remark}

\begin{acknowledgement}
Prashant J. Makadiya acknowledges the support of Government of Gujarat thorugh the SHODH (ScHeme Of Developing High Quality Research) fellowship. Manoj Choudhuri thanks L. Singhal for helpful discussions.
\end{acknowledgement}

\smallskip
\bibliography{ref}
\bibliographystyle{plain}
\end{document}